\documentclass{article}

\usepackage{amsfonts}
\usepackage{amssymb}
\usepackage{amsmath}
\usepackage{enumerate}

\newtheorem{ttt}{Theorem}[section]
\newtheorem{llll}[ttt]{Lemma}
\newtheorem{ccc}[ttt]{Claim}
\newtheorem{eee}[ttt]{Example}
\newtheorem{rrr}[ttt]{Remark}
\newtheorem{sss}[ttt]{Statement}
\newtheorem{ddd}[ttt]{Definition}
\newtheorem{qqq}[ttt]{Question}
\newtheorem{cccc}[ttt]{Corollary}

\newcommand{\bt}{\begin{ttt}}
\newcommand{\bl}{\begin{llll}}
\newcommand{\bc}{\begin{ccc}}
\newcommand{\bex}{\begin{eee}}
\newcommand{\br}{\begin{rrr}\upshape}
\newcommand{\bs}{\begin{sss}}
\newcommand{\bd}{\begin{ddd}\upshape}
\newcommand{\bq}{\begin{qqq}}
\newcommand{\bcor}{\begin{cccc}}

\newcommand{\bp}{\noindent\textbf{Proof. }}

\newcommand{\et}{\end{ttt}}
\newcommand{\el}{\end{llll}}
\newcommand{\ec}{\end{ccc}}
\newcommand{\eex}{\end{eee}}
\newcommand{\er}{\end{rrr}}
\newcommand{\es}{\end{sss}}
\newcommand{\ed}{\end{ddd}}
\newcommand{\eq}{\end{qqq}}
\newcommand{\ecor}{\end{cccc}}

\newcommand{\ep}{\hspace{\stretch{1}}$\square$\medskip}

\newcommand{\lab}[1]{\label{#1}}

\newcommand{\DD}{\mathbb{D}}
\newcommand{\NN}{\mathbb{N}}
\newcommand{\QQ}{\mathbb{Q}}
\newcommand{\RR}{\mathbb{R}}
\newcommand{\ZZ}{\mathbb{Z}}

\newcommand{\al}{\alpha}

\newcommand{\e}{\varepsilon} 
\newcommand{\om}{\omega}

\newcommand{\iH}{\mathcal{H}}
\newcommand{\iI}{\mathcal{I}}

\newcommand{\iP}{\mathcal{P}}
\newcommand{\iS}{\mathcal{S}}

\newcommand{\cat}{^{\mathord{\frown}}}

\begin{document}

\title{Borel sets which are null or non-$\sigma$-finite
for every translation invariant measure}


\author{M\'arton Elekes\thanks{Corresponding author. Partially supported by Hungarian Scientific Foundation grant no.~37758 and F 43620.}\medskip\\
R\'enyi Alfr\'ed Institute\\ 
Budapest, POB 127, H-1364 Hungary\\
emarci@renyi.hu\\
\medskip
and\\ 
\medskip
Tam\'as Keleti\thanks{Partially supported by Hungarian Scientific
Foundation grant no.~F 43620.}\\
E\"otv\"os Lor\'and University\\
Department of Analysis\\ 
Budapest, H-1117 Hungary\\
elek@cs.elte.hu}


\maketitle

\begin{abstract}
We show that the set of Liouville numbers is either null or
non-$\sigma$-finite with respect to every translation invariant Borel
measure on $\RR$, in particular, with respect to every Hausdorff
measure $\iH^g$ with gauge function $g$. This answers a question of
D.~Mauldin. We also show that some other simply defined Borel sets
like non-normal or some Besicovitch-Eggleston numbers, as well as
all Borel subgroups of $\RR$ that are not $F_\sigma$ possess the
above property. We prove that, apart from some trivial cases, the Borel
class, Hausdorff or packing dimension of a Borel set with no such 
measure on it can be arbitrary.
\end{abstract}

\insert\footins{\footnotesize{MSC codes: Primary 28C10; Secondary 28A78,
43A05}}
\insert\footins{\footnotesize{Key Words: translation invariant Borel
measure, Liouville numbers, $\sigma$-finite}} 

\bigskip
\bigskip

\section*{Introduction}
\bigskip

In many branches of mathematics a standard tool is that `nicely
defined' sets admit natural probability measures. For example, limit
sets in the theory of Iterated Function Systems or Conformal Dynamics
as well as self-similar sets in Geometric Measure Theory are usually
naturally equipped with an invariant Borel measure, very often with a
Hausdorff or packing measure. In many situations the sets in
consideration are unbounded, for example periodic, so we cannot hope
for an invariant probability measure. Similarly, the trajectories of
the Brownian motion are of positive $\sigma$-finite $\iH^g$-measure
with probability 1, where the gauge function $g$ is $g(t)=t^2
\log\log\frac{1}{t}$ in case of planar Brownian motion and $g(t)=t^2
\log\frac{1}{t}\log\log\log\frac{1}{t}$ in dimension 3 and
higher. Therefore the natural notion to work with is that of an
invariant positive and $\sigma$-finite measure.

It is natural to ask if there is some sort of unified theory behind
the existence of these measures, for example, one is tempted to ask if
every Borel subset of $\RR^n$ of some `regular structure' is positive
and $\sigma$-finite for some Hausdorff measure, or at least admits a
positive and $\sigma$-finite invariant measure. In
particular, R.~D.~Mauldin (\cite{M2}, \cite{Cs} and see \cite{BD} and
\cite{Bl} for partial and related results) formulated
this question about a specific well-known set of very nice structure;
the set of Liouville numbers, denoted by $L$:
 
\bd
\[
L=\left\{ x\in\RR\setminus\QQ : \forall n\in \NN \ \exists p,q\in \NN
\ (q\geq 2)\textrm{ such that } \left| x-\frac{p}{q} \right|
<\frac{1}{q^n} \right\}.
\]
\ed

\bq
(Mauldin) Is there a translation invariant Borel measure on $\RR$ such
that the set of Liouville numbers is of positive and $\sigma$-finite
measure?
\eq

Note that we of course do not require that the measure be
$\sigma$-finite on $\RR$. Not only because Hausdorff measures are
non-$\sigma$-finite on $\RR$, but also because it is well-known that
every $\sigma$-finite translation invariant Borel measure on the real line is
a constant multiple of Lebesgue measure.

As we will answer this question in the negative, we introduce a
definition.

\bd
A nonempty Borel set $B\subset\RR$ is said to be \emph{immeasurable} if it is
either null or non-$\sigma$-finite for every translation invariant Borel
measure on $\RR$.
\ed

The main result of this paper will be Theorem \ref{Liou} stating that
the set of Liouville numbers is immeasurable. Then in the second part
of the paper we show using various methods that there are other
well-known `nice' immeasurable sets. Specifically, we show that the
set of non-normal numbers, the complement of the set of
Besicovitch-Eggleston numbers, $BE(1,0)$ (one of the
Besicovitch-Eggleston classes itself) are all immeasurable. One of the
main tool is Theorem \ref{group} stating that every Borel but not
$F_\sigma$ additive subgroup of $\RR$ is immeasurable. Using this we
also show that there are immeasurable sets of arbitrary Borel class
(except of course open, as sets of positive Lebesgue measure are
obviously not immeasurable). Similarly, we provide examples of
immeasurable sets of arbitrary Hausdorff or packing dimension.

We note here that it is not only the regular structure of the sets
considered here that makes it difficult to prove immeasurability. Even
it is highly nontrivial to construct some immeasurable set. The two
papers \cite{La} and \cite{Da} containing the two known examples are
entirely devoted to the constructions of the two immeasurable sets.

We would like to point out that in this paper a \emph{Borel measure}
is a measure defined on a $\sigma$-algebra \emph{containing} the Borel
sets. A \emph{group} is always an additive subgroup of
$\RR$. Otherwise the notation we follow throughout the paper can be
found for example in \cite{Ma} or \cite{Fa} and \cite{Ke}. $\lambda$
denotes Lebesgue measure, $\textup{int}(A)$ is the interior of the set
$A$, $A+B=\{a+b:a\in A,b\in B\}$, $A+t=\{a+t:a\in A\}$,
$\textup{card}\left(X\right)$ is the cardinality of the set $X$, a
$G_\delta$ set is the countable intersection of open sets, a
$G_{\delta\sigma}$ set is the countable union of $G_\delta$ sets,
etc. By \emph{Cantor set} we mean a set homeomorphic to the classical
`middle-thirds' Cantor set.

\section{Liouville numbers}

In this section we answer the question of Mauldin.

\bt\lab{Liou}
The set of Liouville numbers is immeasurable; that is, either null or
non-$\sigma$-finite for every translation invariant Borel measure on $\RR$.
\et

It is well-known and easy to check that $L$ is of Lebesgue measure zero, dense,
$G_\delta$ and periodic mod $\QQ$ (that is, $L+q=L$ for every $q\in\QQ$). Hence
the above theorem is a corollary to the following.

\bt\lab{main}
Let $B\subset \RR$ be a nonempty $G_\delta$ set of Lebesgue measure zero. Assume
that $\{t\in\RR:B+t\subset B\}$ is dense in $\RR$. Then $B$ is immeasurable.
\et

The proof of Theorem \ref{main} will be based on two lemmas. The first one 
is reminiscent of \cite{Pr}, where similar results are proved for finite
measures using more complicated methods.

\bl\lab{C}
Let $B$ be a Borel set of Lebesgue measure zero and $\mu$ a Borel measure on
$\RR$ for which $B$ is positive and $\sigma$-finite. Then
\begin{enumerate}[(i)]
\item\lab{fub} $\mu(B\cap(B+t))=0$ for $\lambda$-a.e.~$t$,
\item\lab{int} there exists a Borel set $B'\subset B$ with $\mu(B')>0$ and
$\textup{int}(B'-B')=\emptyset$,
\item\lab{comp} there exists a compact set $C\subset B$ with $\mu(C)>0$ and
$\textup{int}(C-C)=\emptyset$ 
\end{enumerate}
\el

\bl\lab{P} 
Let $B$ be a dense $G_\delta$ set such that $\{t\in\RR:B+t\subset B\}$
is dense in $\RR$, and $C\subset B$ be a compact set with
$\textup{int}(C-C)=\emptyset$. Then there are uncountably many (in fact,
continuum many) disjoint translates of $C$ inside $B$.
\el

It is easy to see that applying Lemma \ref{P} to $C$ of Lemma \ref{C}
\textit{(\ref{comp})} yields Theorem \ref{main} 

In the rest of the section we prove the two lemmas.


\bp 
(Lemma \ref{C}) \textit{(\ref{fub})} Let $\mu:\iS\to[0,\infty]$ be the 
given measure, where $\iS$ is a $\sigma$-algebra of subsets of $\RR$
containing all Borel sets. Define a new measure $\mu_B$ by
\[
\mu_B(S)=\mu(B\cap S) \textrm{ for every } S\in\iS.
\]
$\mu_B$ is clearly a $\sigma$-finite Borel measure on $\RR$. Define
\[
\widetilde{B}=\{(x,y)\in\RR^2:x+y\in B\}.
\]
This is clearly a Borel set, hence $\mu_B\times\lambda$-measurable. As
both measures are $\sigma$-finite, we can apply the Fubini theorem to
$\widetilde{B}$. Vertical sections of $\widetilde{B}$ are of the form
\[
\{y\in\RR:x+y\in B \} = \{y\in\RR:y\in B-x \} = B-x,
\]
therefore are all of Lebesgue measure zero. By Fubini,
[$\mu_B\times\lambda](\widetilde{B})=0$, and so
$\lambda$-a.e.~horizontal section of $\widetilde{B}$ is of
$\mu_B$-measure zero. A horizontal section is $\{x\in\RR:x+y\in B \} =
B-y$, therefore for $\lambda$-a.e.~$y$ we obtain $0=\mu_B(B-y) =
\mu(B\cap (B-y))$. Replacing $y$ by $-t$ yields the result.

\bigskip

\textit{(\ref{int})} By \textit{(\ref{fub})} we can choose a countable dense
set $D\subset\RR$ such that $\mu(B\cap(B+d))=0$ for every $d\in
D$. Define
\[
B'= B \setminus \bigcup_{d\in D} (B+d).
\]
It is easy to check that $\mu(B')=\mu(B)>0$ and $D\cap
(B'-B')=\emptyset$, so the proof of \textit{(\ref{int})} is complete.

\bigskip

\textit{(\ref{comp})} It is sufficient to find a compact set 
$C\subset B'$ of positive $\mu$-measure. Since $B'\subset B$, $B'$ 
is $\sigma$-finite for $\mu$. Let $B'=\cup_{n=0}^\infty S_n$, 
where $S_n\in\iS$ and
$\mu(S_n)<\infty$. Since $\mu(B')>0$, there exists $S=S_n\subset B'$
such that $0<\mu(S)<\infty$. Define
\[
\mu_S(A)=\mu(S\cap A) \textrm{ for every Borel set } A\subset \RR.
\]
Note that in contrast with the above $\mu_B$, this time the measure is defined
only for Borel sets.

$\mu_S$ is clearly a finite measure on the Borel sets, hence inner
regular w.r.t.~compact sets \cite[Thm 17.11]{Ke}. Apply this to $B'$, a
Borel set with $\mu_S(B')=\mu(S\cap B')=\mu(S)>0$, and obtain a
compact set $C\subset B'$ with $\mu_S(C)=\mu(S\cap C)>0$. So
$\mu(C)>0$ follows, and the proof of Lemma \ref{C} is complete.
\ep


\bp 
(Lemma \ref{P}) Let
\[
T=\{t\in\RR:C+t\subset B\}.
\]
$B$ is $G_\delta$, so there are open sets $U_n$ such that
$B=\cap_{n=0}^\infty U_n$. Clearly $C+t\subset \cap_{n=0}^\infty U_n$
iff $C+t\subset U_n$ holds for every $n\in\NN$.  Therefore
$T=\cap_{n=0}^\infty G_n$, where $G_n=\{t\in\RR:C+t\subset U_n\}$. As
$C$ is compact and $U_n$ is open, $G_n$ is open. Note that $G_n$ is
also dense by our assumption.

It is clearly sufficient to construct a Cantor set $P\subset T$ with
the property that $(C+p_0)\cap(C+p_1)=\emptyset$ holds for every pair
of distinct points $p_0,p_1\in P$. We define $P$ via a usual Cantor
scheme as follows. Let $2^n$ stand for the set of 0--1 sequences of length
$n$. We define nondegenerate compact intervals $I_s$ for every
$n\in\NN$ and $s\in 2^n$ by induction on $n$ (we also make sure that
at level $n$ all intervals are of length at most $\frac{1}{n}$). Fix
$I_\emptyset$ such that $I_\emptyset\subset G_0$. Once $I_s$ is
already defined for some $s\in 2^n$, we pick $x\in I_s\cap
G_{n+1}$. As $G_{n+1}$ is open dense and $C-C$ is nowhere dense we can
find $y\in \left[I_s\cap G_{n+1}\right] \setminus
\left[(C-C)+x\right]$. This ensures $(C+x)\cap(C+y)=\emptyset$, as
otherwise $c_0+x=c_1+y$ for some $c_0, c_1 \in C$, so $y=(c_0-c_1)+x$,
a contradiction. By compactness we can find disjoint $I_{s\cat 0},
I_{s\cat 1} \subset I_s\cap G_{n+1}$ such that $x\in I_{s\cat 0}$,
$y\in I_{s\cat 1}$ and $(C+I_{s\cat 0}) \cap (C+I_{s\cat 1})
=\emptyset$.
 
Now define
\[
P=\bigcap_{n=0}^\infty \bigcup_{s\in 2^n} I_s.
\]
Then clearly $P$ is a Cantor set, $P\subset T=\cap_{n=0}^\infty G_n$
as $I_s\subset G_n$ for every $n$ and $s\in 2^n$. Moreover, if $p_0,
p_1\in P$, $p_0\neq p_1$ then there are $n$ and $s\in 2^n$ such that
$p_0\in I_{s\cat 0}$ and $p_1\in I_{s\cat 1}$ (or the other way
around), but then $(C+p_0)\cap(C+p_1)=\emptyset$ holds since
$(C+I_{s\cat 0}) \cap (C+I_{s\cat 1}) =\emptyset$. This completes the
proof of Lemma \ref{P}.
\ep

\br
Theorem \ref{main} holds (with essentially the same proof) in
$\RR^n$, more generally in locally compact abelian Polish groups as
well. Of course, Lebesgue measure is replaced by Haar measure.

No assumptions of Theorem \ref{main} can be omitted. The empty set is
not immeasurable by definition. The example of $\QQ$ and the counting
measure shows that the $G_\delta$ assumption is essential. Sets of
positive Lebesgue measure are clearly not immeasurable. The density of
$\{t\in\RR:B+t\subset B\}$ is also important, as the example $B=L\cup
C_{1/3}$ shows, where $C_{1/3}$ is the classical middle-thirds Cantor
set. Indeed, $\log 2/\log 3$-dimensional Hausdorff measure is positive
and finite on $B$, since it is easy to see that $\dim_H L=0$; that is,
Hausdorff dimension of $L$ is zero (see \cite{Ma} for the definition).

Theorem \ref{main} provides a lot of immeasurable sets in the
following sense.  Let $A$ be a (nonempty) set of Lebesgue measure zero, and let
$B$ be a $G_\delta$ set of measure zero containing $A+\QQ$. Then
$\cap_{q\in\QQ} (B+q)$ fulfills all conditions of Theorem \ref{main},
hence it is immeasurable. Clearly $A\subset \cap_{q\in\QQ} (B+q)$,
therefore every Lebesgue nullset can be covered by an immeasurable
set. Combining this with Lemma \ref{triv} \textit{(\ref{union})} we obtain that
the $\sigma$-ideal generated by the immeasurable sets is the Lebesgue
null ideal.
\er

\section{Non-normal numbers, classes of
Besicovitch-Eggleston numbers and Borel groups}

In this section we provide some more examples of well-known sets that are
immeasurable. 

The following definition is classical.

\bd
A real number $x\in [0,1)$ is called \emph{normal} if its decimal expansion
$x=0.d_1 d_2 \ldots$ satisfies
\[
\lim_{n\to\infty} \frac{\sum_{i=1}^n d_i}{n} = 4.5.
\]
A real number $x\in\RR$ is called normal if its fractional part, $\{x\}$ is
normal. Note that for negative numbers this is \emph{not} the usual decimal
expansion.
The set of normal numbers is denoted by $N$.
\ed

It is well-known and easy to check using the Strong Law of Large Numbers that
$\lambda$-a.e.~real number is normal. 

Our next goal is to prove that the set of non-normal numbers,
$\RR\setminus N$ is immeasurable. Note that the Hausdorff dimension of
this set is 1 (see \cite{BS}), while we already mentioned that $\dim_H
L=0$.  Moreover, $\RR\setminus N$ is not $G_\delta$, since its exact
Borel class is $G_{\delta\sigma}$ (that is, $\RR\setminus N \in
G_{\delta\sigma} \setminus F_{\sigma\delta}$, see \cite[Ex 23.7]{Ke}).

Before the proof we need one more easy lemma.

\bl\lab{triv}
\begin{enumerate}[(i)]
\item\lab{union} The class of immeasurable sets is closed under countable
unions.
\item\lab{count} If the symmetric difference $A \triangle B$ is countable then
$A$ is immeasurable iff $B$ is immeasurable.
\end{enumerate} 
\el

\bp
\textit{(\ref{union})} This clearly follows from the definition.

\textit{(\ref{count})} If $A$ and $B$ are nonempty countable then 
counting measure works for both. 
Otherwise all measures we consider are continuous; that is, singletons
are of measure zero.
\ep

\bt\lab{norm}
The set $\RR\setminus N$ of non-normal numbers is immeasurable; that is,
either null or non-$\sigma$-finite for every translation
invariant Borel measure on $\RR$.
\et

\bp
By Lemma \ref{triv} \textit{(\ref{count})} it is enough to show that $(\RR\setminus
N)\setminus \QQ$ is immeasurable. 
Write this set as
\[
(\RR\setminus N)\setminus \QQ = \bigcup_{k=1}^\infty A_k \cup
\bigcup_{k=1}^\infty B_k,
\]
where $A_k = \left\{x\in\RR\setminus\QQ: \limsup_{n\to\infty}
\frac{\sum_{i=1}^n d_i}{n} \geq 4.5+\frac{1}{k} \right\}$ and similarly $B_k =
\left\{x\in\RR\setminus\QQ: \liminf_{n\to\infty} \frac{\sum_{i=1}^n d_i}{n}
\leq 4.5-\frac{1}{k} \right\}$. By Lemma \ref{triv} (\ref{union}) it is enough
to prove that the sets $A_k$ and $B_k$ are immeasurable for every $k$. The
other case being similar we check this only for $A_k$. 

We show that Theorem \ref{main} applies. As a.e. number is normal,
$A_k$ is of Lebesgue measure zero. Moreover, it is easy to see that it
is dense and periodic modulo a dense set, namely the set of real
numbers with finitely many nonzero digits. Therefore what remains to
be checked is that $A_k$ is $G_\delta$.
\[
A_k = \bigcap_{l=1}^\infty \bigcap_{m=1}^\infty \bigcup_{n=m}^\infty
\left\{x\in\RR\setminus\QQ: \frac{\sum_{i=1}^n d_i}{n} >
4.5+\frac{1}{k}-\frac{1}{l} \right\},
\]
and $\left\{x\in\RR\setminus\QQ: \frac{\sum_{i=1}^n d_i}{n} >
4.5+\frac{1}{k}-\frac{1}{l} \right\}$ is clearly an open subset relative to
$\RR\setminus\QQ$, hence $A_k$ is $G_\delta$ relative to
$\RR\setminus\QQ$, hence $A_k$ is $G_\delta$ in $\RR$ as well.
\ep

Our next example concerns the so called Besicovitch-Eggleston numbers. These
are the numbers for which `the asymptotic frequency of every digit exists':
Let $\{x\}=0.d_1 d_2 \ldots$ be the expansion of the fractional part of the
real number $x$ in base $b\geq 2$. 
If there are two such expansions, we choose the one with
finitely many nonzero digits.

\bd
Let $\alpha_0,\ldots,\alpha_{b-1}$ be real numbers such that $\alpha_d\geq 0$
for every $d=0,\ldots,b-1$ and $\sum_{d=0}^{b-1} \alpha_d =1$. 
\[
BE(\alpha_0,\ldots,\alpha_{b-1}) =
\]
\[ \left\{
x\in\RR:
\lim_{n\to\infty}
\frac{\textup{card}\left(\{i:1\leq i\leq n, d_i=d \} \right)}{n} = \alpha_d
\textrm{ for every } d=0,\ldots,b-1 
\right\}.
\]
The union of all these classes (for a fixed positive integer $b$) is denoted
by $BE_b$.
\ed

\bt
For every integer $b\geq 2$ the set $\RR\setminus BE_b$ is immeasurable.
\et

\bp
The proof is very similar to that of Theorem \ref{norm}. We can again restrict
ourselves to the set of
irrational numbers. Write 
\[
(\RR\setminus BE_b)\setminus\QQ = \bigcup_{d=0}^{b-1}
\bigcup_{\substack{p<q \\ p,q\in \QQ\cap [0,1]}} \left( C_q^d\cap
D_p^d\right),
\]
where 
\[
C_q^d = \left\{x\in\RR\setminus\QQ:
\limsup_{n\to\infty} \frac{\textup{card}\left(\{i:1\leq i\leq n, d_i=d \}
\right)}{n} \geq q \right\}
\]
and 
\[
D_p^d = \left\{x\in\RR\setminus\QQ:
\liminf_{n\to\infty} 
\frac{\textup{card}\left(\{i:1\leq i\leq n, d_i=d \} \right)}{n} \leq p
\right\}.
\]
It is again easy to
see using the Strong Law of Large Numbers that Lebesgue a.e.~number is in
$BE_b$, namely in $BE(\frac{1}{b},\ldots,\frac{1}{b})$, hence every $C_q^d$ and
$D_p^d$ is of Lebesgue
measure zero. Moreover, these sets are invariant under all translations 
from the dense set
of numbers with finitely many nonzero digits. So it is enough to show
that $C_q^d$ and $D_p^d$ are $G_\delta$ in $\RR\setminus\QQ$. But this
is clear, as
\[
C_q^d = \bigcap_{l=1}^\infty \bigcap_{m=1}^\infty \bigcup_{n=m}^\infty
\left\{x\in\RR\setminus\QQ: 
\frac{\textup{card}\left(\{i:1\leq i\leq n, d_i=d
\}\right)}{n} > q-\frac{1}{l} 
\right\},
\]
and similarly for $D_p^d$.
\ep

So far all immeasurable sets in this paper were comeager. Our next
example is easily seen to be meager. Not surprisingly we need new
methods to prove immeasurability here.

\bd
Let $A$ be a set of natural numbers. The expression
\[
\lim_{n\to\infty}
\frac{\textup{card}\left(A\cap[1,n] \right)}{n}
\]
(if exists) is called the \emph{density of $A$}.
\ed

It is easy to see that sets of zero density form an ideal on the natural
numbers. Dually, sets of density 1 form a filter.

\bt\lab{be}
$BE(1,0)$ is immeasurable; that is, the set of real numbers with the property
that in the dyadic expansion the 1's form a set of zero density is either null
or non-$\sigma$-finite for every translation invariant Borel measure on $\RR$.
\et

This result will follow from the following more general statement.

\bt\lab{gengroup}
Let $B$ be a Borel set such that $B+B\subset B$. If $B-B$, the
group generated by $B$, is not $F_\sigma$ then $B$ is immeasurable.
\et


\bp
Assume towards a contradiction that $\mu$ is a translation invariant
measure for which $B$ is positive and $\sigma$-finite.  As in the
proof of Lemma \ref{C} \textit{(\ref{comp})} we obtain that there exists a
compact set $C\subset B$ of positive $\mu$-measure. We get a
contradiction by showing that $C$ has uncountably many disjoint
translates in $B$.  We define the transfinite sequence of numbers
$\{t_\al:\al<\om_1\}$ by transfinite induction, so that the sets
$\{C+t_\al:\al<\om_1\}$ are pairwise disjoint. Clearly $(C+t_\al) \cap
(C+t_\beta) = \emptyset$ iff $t_\al \notin C-C+t_\beta$. Throughout
the induction we choose the numbers from $B$, which implies $C+t_\al
\subset B$ by our assumption, so what remains to show is that the
induction cannot get stuck at some countable ordinal. At step $\al$
our task is to find $t_\al\in B$ so that $t_\al \notin C-C+t_\beta$
for every $\beta<\al$. These latter sets are in $B-B$ because this is
a group, and we claim that they cannot cover $B$. Indeed, if $B\subset
\cup_{\beta<\al} (C-C+t_\beta)$ then $B-B=\cup_{\beta,\gamma<\al}
[(C-C+t_\beta)-(C-C+t_\gamma)]$, so $B-B$ is $F_\sigma$, which is 
a contradiction.
\ep

Before we finish the proof of Theorem \ref{be} we need one more lemma, which 
provides an example of an immeasurable group.

\bl\lab{G}
Denote 
\[
G=\left\{x\in\RR:d_n=d_{n+1} \textrm{ for
all } n \textrm{ but a set of density zero}\right\},
\]
where $\{x\}=0.d_1 d_2\ldots$ is the dyadic expansion of the fractional part 
of $x$. Then
\begin{enumerate}[(i)]
\item\lab{G1} $BE(1,0)-BE(1,0)=G$,
\item\lab{G2} $G$ is a Borel but not $F_\sigma$ group.
\end{enumerate} 
\el

\bp
\textit{(\ref{G1})}
Put $B=BE(1,0)$ and define $S(x)=\{n\in\NN:d_n\neq d_{n+1}\}$, where 
$\{x\}=0.d_1 d_2\ldots$ is the dyadic expansion of the fractional part 
of $x$. Then 
$G=\{x\in\RR: S(x)\textrm{ is of density zero}\}$.
We have to show that $B-B = G$. Note that all
dyadic rationals (numbers with only finitely many nonzero digits) are in $B$
as we chose the expansion that is eventually zero. 

First we check $B-B\subset
G$. Fix $b_1,b_2\in B$, and write $b_1-b_2=b_1+(-b_2)$. 
If $b_2$ is a dyadic
rational then so is $-b_2$, hence all but finitely many digits of $b_1+(-b_2)$
and $b_1$ coincide, hence $b_1-b_2\in B$. So we can assume that $b_2$ is not a
dyadic rational, and then we obtain the expansion of (the fractional part of)
$-b_2$ by replacing 0's with 1's and vice versa. Let $d_i$, $e_i$ and $f_i$  
be the
digits of the expansions of $b_1$, $-b_2$ and $b_1+(-b_2)$, respectively. 
Then the set
$A=\{i: d_i=d_{i+1}=0, e_i=e_{i+1}=1\}$ is of density 1. 
If $i\in A$ then
either $f_i=f_{i+1}=0$ or $f_i=f_{i+1}=1$ (depending on whether there is 
a carried digit from the $i+2^{nd}$ place), hence $i\notin S(b_1+(-b_2))$.
Therefore $S(b_1+(-b_2))$ is of density zero
and so $b_1+(-b_2) \in G$.


Now we show that $G\subset B-B$. Given $x\in G$ with expansion $\{x\}=0.d_1
d_2\ldots$ we construct $b\in B$ with expansion $\{b\}=0.e_1 e_2\ldots$ such
that $x+b\in B$. We can assume that $x$ is not a dyadic rational, otherwise
$b=0$ works. Define $b$ as follows: let $e_i=1$ iff $d_i=1$ and $d_{i+1}=0$.
Since $x\in G$ we have that $S(x)$ is of density zero, therefore $b\in
B$. Let $\{x+b\}=0.f_1 f_2\ldots$. It suffices to check that $f_i=1$ 
iff $d_i=0$ and $d_{i+1}=1$, hence then $\{i:f_i=1\}$ is a subset of $S(x)$, 
therefore is of density zero, thus $x+b\in B$. Suppose first that 
$d_{i+1}=0$. Then $e_{i+1}=0$, so there is no carried digit from this place. 
Moreover $e_i=d_i$, so $f_i=0$. Suppose now that $d_{i+1}=1$. Let $k$ be 
the minimal integer such that $i+2\leq k$ and $d_k=0$. Then $e_j=0$ for 
$i\leq j \leq k-2$ and for $j=k$, and $e_{k-1}=1$. Consequently, there is 
no carried digit from the $k^{th}$ place, so $f_{k-1}=0$ and there is a 
carried digit from the $k-1^{st}$ place, and so on there is a carried 
digit for every $i+1 \leq j \leq k-1$, and at the end $f_i=1$ iff $d_i=0$.

\textit{(\ref{G2})}
$B=BE(1,0)$ and $G$ are clearly Borel sets. As $B$ is easily seen to 
be closed under addition, $G=B-B$ is a group. So we have to check 
that $G$ is not $F_\sigma$. 

Suppose towards a contradiction that it is. First we 
construct a continuous map
$f:[0,1)\setminus\DD\to[0,1)\setminus\DD$ (where $\DD$ is the set of
dyadic rationals), with the property that
$f^{-1}(G\cap([0,1)\setminus\DD))=B\cap([0,1)\setminus\DD)$. For
$x=0.d_1 d_2\ldots$ construct the digits of $f(x)=0.e_1 e_2\ldots$ as
follows. Let $d_1=e_1$, and then define $e_n$ recursively such that
$e_{n+1}=e_n$ iff $d_n=0$. It is easy to check that $f$ satisfies the
requirements. By assumption $G$ is $F_\sigma$, hence
$G\cap([0,1)\setminus\DD)$ is $F_\sigma$ in $[0,1)\setminus\DD$, so
$B\cap([0,1)\setminus\DD)$ is also $F_\sigma$ in
$[0,1)\setminus\DD$. Then $B\cap [0,1)$
is $G_{\delta\sigma}$, which is a contradiction, as by \cite[Ex 23.7]{Ke} 
it is not $G_{\delta\sigma}$.
\ep

\bp
(Theorem \ref{be}) As we already mentioned, $BE(1,0)$ is clearly Borel and 
closed under 
addition, hence Lemma \ref{G} shows that Theorem \ref{gengroup} applies.
\ep

We also formulate the following immediate corollary to Theorem 
\ref{gengroup}, which 
is interesting in its own right.

\bt\lab{group}
If $G\subset\RR$ is a Borel but not $F_\sigma$ additive subgroup then $G$ is immeasurable. $\square$
\et

\br
In fact this theorem holds in every Polish group, and $F_\sigma$ can
be weakened to $K_\sigma$ (that is, countable union of compact sets). 
In case the group is not abelian we have to
use left-translations everywhere (including the definition of
immeasurable.)
\er

\section{Immeasurable sets with arbitrary Borel class, Hausdorff or packing 
dimension}

A nonempty open set is of positive Lebesgue measure, hence cannot be
immeasurable. Our next theorem shows that we can find immeasurable
sets in all other Borel classes.

For $1\leq\al<\omega_1$ denote $\Sigma_\alpha^0$, $\Pi_\alpha^0$ and
$\Delta_\alpha^0$ the additive, multiplicative and ambiguous Borel
classes, respectively (see e.g.~\cite{Ke}). We say that a Borel set
$B$ is of exact class $\Sigma_\alpha^0$\ ($\Pi_\alpha^0$) if $B\in
\Sigma_\alpha^0 \setminus \Pi_\alpha^0 \ (B\in
\Pi_\alpha^0\setminus\Sigma_\alpha^0)$. We say that $B$ is of exact
class $\Delta_\alpha^0$ if $B\in\Delta_\alpha^0$ but
$B\notin\Sigma_\beta^0 \cup \Pi_\beta^0$ for $\beta<\al$.

\bt 
For every $1\leq\al<\omega_1$ there exist immeasurable sets of exact Borel 
class $\Pi_\alpha^0$. There exist immeasurable sets of exact Borel 
class $\Sigma_\alpha^0$ and $\Delta_\alpha^0$ iff $2\leq\al<\omega_1$.
\et

\bp
$\Delta_1^0$: Impossible, as open sets are not immeasurable.

$\Sigma_1^0$: Impossible, for the same reason.

$\Pi_1^0$: Davies \cite{Da} constructed an immeasurable Cantor set $D$.

$\Delta_2^0$: $D$ minus a point is immeasurable by Lemma \ref{triv} 
\textit{(\ref{count})}.

$\Sigma_2^0$: $D+\QQ$ is immeasurable by Lemma \ref{triv}
(\ref{union}).  To see that its exact class is $F_\sigma$, first note that as
$D$ is immeasurable, it is of Lebesgue measure zero. Moreover it is
closed, hence nowhere dense, therefore $D+\QQ$ is a dense meager
$F_\sigma$ set. Thus it cannot be $G_\delta$, as dense $G_\delta$ sets
are comeager.

$\Pi_2^0$: The set of Liouville numbers is immeasurable, see Theorem
\ref{main}. Larman's example \cite{La} is also $G_\delta$. Both
examples are easily seen to be non-$F_\sigma$ using the fact that
nonmeager $F_\sigma$ sets have interior points.

All other classes: These examples can be chosen to be groups: in
\cite{M1} Borel groups of exact class $\Pi_\alpha^0$,
$\Sigma_\alpha^0$ and $\Delta_\alpha^0$ are constructed (for
$3\leq\al<\omega_1$), and these are all immeasurable by Theorem
\ref{group}.
\ep

Our next theorem shows that there are lots of immeasurable
sets from the viewpoint of Hausdorff dimension as well. 


\bt\lab{Haus}
For every $0\leq\al\leq 1$ there exists an immeasurable subset of $\RR$ 
of Hausdorff dimension $\al$. In particular, these sets can be 
chosen to be additive subgroups of $\RR$.
\et

\bp
First we consider the case $\al\neq 1$.

In \cite{EV} groups $G_\al\subset\RR$ of Hausdorff dimension $\al$ are
constructed for every $0<\al<1$. In fact, the proof also yields a
nontrivial group $G_0$, but for our purposes it is sufficient to put
$G_0=\{0\}$. We will also use the fact that
$G_\al\subset G_\beta$ for $\al<\beta$. 

Unfortunately $G_\al$ is $F_\sigma$. Our goal is to find a group $H$
`sufficiently orthogonal to $G_\al$' such that $G'_\al=G_\al+H$ is a
non-$F_\sigma$ group of dimension $\al$.



The following lemma is probably well-known, however, we were unable to
find suitable references for its second half, so we include a proof
here.

\bl\lab{typic}
Let $M$ be a meager subset of $\RR\setminus\{0\}$. Then the typical
compact subset $C$ of $\RR$ (that is, comeager many elements of the
space of compact subset of $\RR$ endowed with the Hausdorff metric)
possesses the following properties.
\begin{enumerate}[(i)]
\item\lab{indep} $C$ is a Cantor set which is linearly independent
over the rationals.
\item\lab{inters} The group generated by $C$ is disjoint from $M$.
\end{enumerate}
\el

\bp
\textit{(\ref{indep})} This is essentially \cite[Ex 19.2]{Ke} which 
follows from \cite[Thm 19.1]{Ke}.

\textit{(\ref{inters})} The complement of $M$ contains a set that can be
written as the intersection of countably many dense open sets. Hence it is
enough to prove that for a fixed dense open set $U$, the group
generated by $C$ is in $U$.

For $(n_1,\ldots,n_k,n'_1,\ldots,n'_l)\in (\NN\setminus\{0\})^{k+l}$
define
\[
C(n_1,\ldots,n_k;n'_1,\ldots,n'_l) =
\]
\[
=\{n_1 c_1 + \ldots + n_k c_k - n'_1 c'_1 -\ldots - n_l c'_l : c_i,
c'_j\in C \textrm{ are all distinct} \}.
\]
Set $C(\emptyset)=\{0\}$.  The union of all these sets as
$(n_1,\ldots,n_k,n'_1,\ldots,n'_l)$ ranges over the countable sets
$(\NN\setminus\{0\})^{k+l}$ is the group generated by $C$, hence it
suffices to prove that
\[
C(n_1,\ldots,n_k;n'_1,\ldots,n'_l)\subset U
\]
for the typical $C$. Define the map $f:\RR^{k+l}\to\RR$ by
\[
f(x_1,\ldots,x_k,x'_1,\ldots,x'_l)=n_1 x_1 + \ldots + n_k x_k - n'_1
x'_1 -\ldots - n_l x'_l.
\]
This map is clearly continuous and open, therefore $f^{-1}(U)$ is
dense and open in $\RR^{k+l}$. Then \cite[Thm 19.1]{Ke} states precisely
that for the typical $C$ the set
\[
(C)^{k+l}=\{(c_1,\ldots,c_k;c'_1,\ldots,c'_l) : c_i, c'_j\in C
\textrm{ are all distinct} \}
\]
is contained in $f^{-1}(U)$. Therefore
$C(n_1,\ldots,n_k;n'_1,\ldots,n'_l)\subset U$, which finishes the
proof of the lemma.
\ep

The last easy statement we need before we define our group $H$ is that
every Cantor subset $K$ of $\RR$ contains a Cantor set $K_0$ of packing
dimension zero. We need this additional step as the typical compact
set has packing dimension 1. See \cite{Ma} for the definition. Packing
dimension (which is the same as upper packing dimension) of a set $B$
will be denoted by $\dim_p B$, upper Minkowski dimension of $B$ is
denoted by $\overline{\dim}_M B$. We will use the well-known
inequality $\dim_p B\leq\overline{\dim}_M B$.

In order to prove this statement, we perform a Cantor-type
construction as follows. We can clearly assume $K\subset [0,1]$. 
For every integer
$m\geq 1$ divide $[0,1]$ into $2^m$ subintervals of length
$\frac{1}{2^m}$.  One can check that it is possible to choose
recursively for every $m$ a nonempty subfamily $\iI_m$ of these intervals such
that
\begin{itemize}
\item $\cup \iI_m \supset \cup \iI_n$ for $m<n$,
\item $\textup{card}(\iI_m)\leq m$
\item $\forall I\in\iI_m \ \exists n\in\NN \textrm{ and } J, J'\in \iI_n, J\neq J', \textrm{ such that } J, J'\subset I$,
\item $\forall I\in\iI_m \ \textup{int}(I)\cap K\neq \emptyset$.
\end{itemize}
Define $K_0=\bigcap_{m=1}^\infty\cup \iI_m$.
This set is easily seen to be a Cantor set, and using the ``box-counting'' version of the upper Minkowski dimension (see \cite[page 78]{Ma}) yield $\overline{\dim}_M K_0=0$. Hence $\dim_p K_0=0$.

Now we complete the proof of Theorem \ref{Haus}. The only nonmeager
Borel subgroup of $\RR$ is $\RR$ itself (since by \cite[page 93]{Ox} if 
$B$ is a nonmeager Borel set than $B-B$ contains an interval). 
Therefore the sets $G_{1-\frac{1}{n}}
\setminus \{0\}$ are meager. Apply Lemma \ref{typic} to
$M=\cup_{n\in\NN} (G_{1-\frac{1}{n}} \setminus \{0\})$ and obtain a
compact set $C$. By the previous argument we can choose a Cantor set
$C_0\subset C$ such that $\dim_p C_0=0$. Fix a Borel but not
$F_\sigma$ subset $B\subset C_0$, and define $H$ as the group
generated by $B$. We have to check that $G'_\al=G_\al+H$ is a Borel
but not $F_\sigma$ group of Hausdorff dimension $\al$
for every $0\leq\al<1$. 

As $G_\al\subset G_\beta$ for $\al<\beta$, we obtain $G_\al\cap H =\{0\}$
for every $0\leq\al<1$.

First we show that $H$ is Borel. As $H$ is the group generated by $B$,
it is the union of all the sets $B(n_1,\ldots,n_k;n'_1,\ldots,n'_l)$
(the notation was defined in the proof of Lemma \ref{typic}). Therefore it is
sufficient to show that these latter sets are all Borel. But since $B$
is linearly independent, such a set is the continuous one-to-one image
of the Borel set $(B)^{k+l}$ (as above, using the notation of
the proof of Lemma \ref{typic}), hence Borel.

Similarly, the natural map from the Borel set $G_\al\times
H\subset\RR^2$ onto $G'_\al=G_\al+H$ is continuous. Moreover, it is 
one-to-one, since $g_1+h_1=g_2+h_2$ implies $g_1-g_2=h_2-h_1$ and 
$G_\al$ and $H$ are groups satisfying $G_\al\cap H =\{0\}$. Hence 
$G'_\al$ is Borel.

In order to show that $G'_\al$ is not $F_\sigma$ will prove
$G'_\al\cap C_0=B$, which is clearly sufficient as $C_0$ is compact
and $B$ is not $F_\sigma$. Suppose $g'=g+h\in C_0$, where $g\in G_\al,
h\in H$. We want to show that $g'\in B$.  Clearly, $g=g'-h\in
C_0-H$. But $C_0-H$ is a subset of the group generated by $C$, and
$G_\al\subset M$, so by Lemma \ref{typic} \textit{(\ref{inters})} we obtain
$g=0$. Hence $h\in C_0$, and $h$ is in the group generated by
$B\subset C_0$. But $C_0$ is linearly independent, therefore $h\in B$
and we are done, as $h=g'$.


What remains to be shown is that $\dim_H G'_\al =
\al$. Obviously $\al = \dim_H G_\al \leq \dim_H G'_\al$, so it
suffices to prove $\dim_H G'_\al\leq\al$. It is well-known \cite[Thm 8.10]{Ma} that $\dim_p (X\times Y)\leq \dim_p X + \dim_p Y$, which implies
$\dim_p (X_1\times \ldots \times X_n)\leq \dim_p X_1 +\ldots +\dim_p
X_n$.  First we prove $\dim_p H=0$. As $H$ is the countable union of
the sets $B(n_1,\ldots,n_k;n'_1,\ldots,n'_l)$, it is sufficient to show that 
$\dim_p B(n_1,\ldots,n_k;n'_1,\ldots,n'_l)=0$. This set can be
covered by a Lipschitz image of the set $B^{k+l}$, and Lipschitz
images do not increase dimension, so it is enough to check $\dim_p
B^{k+l}=0$. But this clearly holds as $B\subset C_0$, so $\dim_p B\leq \dim_p C_0=0$, hence $\dim_p (B\times \ldots \times
B)\leq \dim_p B +\ldots +\dim_p B = 0+\ldots+0 =0$. So $\dim_p H=0$.

Another part of \cite[Thm 8.10]{Ma} states that $\dim_H
(X\times Y)\leq \dim_H X + \dim_p Y$. Choosing $X=G_\al$ and $Y=H$
yields $\dim_H (G_\al\times H)\leq \al+0=\al$. But $G'_\al$ is a
Lipschitz image of $G_\al\times H$, hence $\dim_H G'_\al\leq\al$. 
This completes the proof for $\al<1$.

Finally, we have to deal with the case $\al=1$. As $G'_{1-\frac{1}{n}}$
is an increasing sequence of immeasurable groups of Hausdorff
dimension $1-\frac{1}{n}$, the union $\cup_{n=1}^\infty G_{1-\frac{1}{n}}$
is a group of Hausdorff dimension $1$ which is immeasurable by Lemma
\ref{triv} (\ref{union}).
\ep

\br
In \cite{EV} only $\dim_H G_\al=\al$ is proved, but $\dim_p G_\al=\al$
also holds (see sketch of proof below). Therefore one can easily check
that a minor modification of the above proof yields $\dim_p
G'_\al=\al$, hence Theorem \ref{Haus} is also valid for packing
instead of Hausdorff dimension.

Now we sketch the proof of $\dim_p G_\al=\al$. Since $\al=\dim_H
G_\al\leq\dim_p G_\al$, it is enough to prove $\dim_p G_\al\leq\al$.
Write $G_\al=\cup_{k_0\in\NN}\cup_{\kappa\in\NN}\cup_{n\in\ZZ}
(G_\al(k_0,\kappa)+n)$, where $G_\al(k_0,\kappa)$ is the set of points
in $[0,1]$ satisfying the requirement in the definition of $G_\al$ for
fixed $k_0$ and $\kappa$ (see the construction in \cite{EV}). It is
sufficient to show $\dim_p G_\al(k_0,\kappa)\leq\al$ for fixed $k_0$
and $\kappa$. This set is a Cantor set, and it is not hard to see that
the `natural measure' on it (we always split the measure into equal
parts in the Cantor-type construction) satisfies all conditions of
\cite[Theorem 6.11 (2)]{Ma} with $A=G_\al(k_0,\kappa)$, $\lambda=1$ and
$s=\al+\e$, where $\e>0$ is arbitrary. Hence $\dim_p
G_\al(k_0,\kappa)\leq\al$.
\er

\section{Concluding remarks}

There are numerous questions arising naturally concerning
immeasurability.  Is it really weaker to require some translation
invariant Borel measure than a Hausdorff or packing measure? The
answer to the second version turns out to be positive, as Peres
\cite{Pe} showed that some of the so called Bedford-McMullen carpets
are of zero or non-$\sigma$-finite packing measure for every gauge
function, while by recent results of the authors
\cite{EK} none of these carpets are immeasurable.  We do not know the
answer to the other version: Does the existence of a translation
invariant Borel measure that is positive and $\sigma$-finite on a
Borel set $B$ imply the existence of a Hausdorff measure with a
suitable gauge function that is positive and $\sigma$-finite on $B$?

We may be overlooking something simple, but the
following questions also seem to be open. We say that a Borel set $B$ is
\emph{measured} if it is not immeasurable; that is, there is a
translation invariant Borel measure for which $B$ is positive $\sigma$-finite.
Is the union of two measured sets also measured? (It can be shown using the construction in \cite{Da} that 
this is false for countably many sets.) Is the set of
Liouville numbers a finite/countable union of measured sets? Is every
Borel set a finite/countable union of measured sets?

\bigskip

\noindent
\textbf{Acknowledgement} Originally, the results of this paper were
proven for measures defined on the Borel sets. The
authors are indebted to D.~H.~Fremlin for pointing out how to
generalize the results for measures defined on $\sigma$-algebras
containing the Borel sets. We are also grateful to R.~D.~Mauldin for
helpful discussions.










\end{document}